\documentclass[11pt]{article}
\usepackage{amsfonts,amsbsy,amsmath,amssymb}
\usepackage{graphics}
\usepackage[lmargin=3.0cm, rmargin=3.0cm, tmargin=3.0cm, bmargin=3.0cm]{geometry}

\numberwithin{equation}{section}

\newcommand{\iso}{\cong} 
\newcommand{\notiso}{\not\cong} 
\newcommand{\size}[1]{\Vert #1 \Vert}

\newcommand{\scaps}[1]{{\scshape #1}}
\newcommand{\bfm}{\textbf}
\newcommand{\mcal}{\mathcal}

\newcommand{\mbs}{\boldsymbol}

\newcommand{\B}{\mcal{B}}
\newcommand{\C}{\mcal{C}}

\newcommand{\F}{\mcal{F}}

\renewcommand{\H}{\mcal{H}}

\renewcommand{\P}{\mcal{P}}

\newcommand{\T}{\mcal{T}}

\newcommand{\bnd}[2]{{\rm bnd}_{#2}#1}

\newcommand{\sm}{\setminus}

\def\proof{\par\noindent{{\sl Proof.\ }}}
\def\QED{$\blacksquare$}
\def\inQED{$\square$}

\newenvironment{statement}
{
\refstepcounter{equation} 
\ \\
\noindent
\begin{it}
\noindent
\emph{\bfm{~\theequation.}}
}
{\end{it}\ \\}

\newcommand{\sref}[1]{\bfm{{\ref{#1}}}}

\makeatletter
\def\section{\@ifstar\unnumberedsection\numberedsection}
\def\numberedsection{\@ifnextchar[
  \numberedsectionwithtwoarguments\numberedsectionwithoneargument}
\def\unnumberedsection{\@ifnextchar[
  \unnumberedsectionwithtwoarguments\unnumberedsectionwithoneargument}
\def\numberedsectionwithoneargument#1{\numberedsectionwithtwoarguments[#1]{#1}}
\def\unnumberedsectionwithoneargument#1{\unnumberedsectionwithtwoarguments[#1]{#1}}
\def\numberedsectionwithtwoarguments[#1]#2{%
  \ifhmode\par\fi
  \removelastskip
  \vskip 3ex\goodbreak
  \refstepcounter{section}%
  \noindent
  \leavevmode
  \begingroup
  \bfseries
  \S \thesection\ 
  #2.\quad
  \endgroup
  \addcontentsline{toc}{section}{%
    \protect\numberline{\thesection}%
    #1}%
  }
\def\unnumberedsectionwithtwoarguments[#1]#2{%
  \ifhmode\par\fi
  \removelastskip
  \vskip 3ex\goodbreak
  \noindent
  \leavevmode
  \begingroup
  \bfseries
  #2.\quad
  \endgroup
  \addcontentsline{toc}{section}{%
    #1}%
  }

\makeatletter
\def\subsection{\@ifstar\unnumberedsubsection\numberedsubsection}
\def\numberedsubsection{\@ifnextchar[
  \numberedsubsectionwithtwoarguments\numberedsubsectionwithoneargument}
\def\unnumberedsubsection{\@ifnextchar[
  \unnumberedsubsectionwithtwoarguments\unnumberedsubsectionwithoneargument}
\def\numberedsubsectionwithoneargument#1{\numberedsubsectionwithtwoarguments[#1]{#1}}
\def\unnumberedsubsectionwithoneargument#1{\unnumberedsubsectionwithtwoarguments[#1]{#1}}
\def\numberedsubsectionwithtwoarguments[#1]#2{%
  \ifhmode\par\fi
  \removelastskip
  \vskip 3ex\goodbreak
  \refstepcounter{subsection}%
  \noindent
  \leavevmode
  \begingroup
  \bfseries
  \S \thesubsection\ 
  #2.\quad
  \endgroup
  \addcontentsline{toc}{subsection}{%
    \protect\numberline{\thesubsection}%
    #1}%
  }
\def\unnumberedsubsectionwithtwoarguments[#1]#2{%
  \ifhmode\par\fi
  \removelastskip
  \vskip 3ex\goodbreak
  \noindent
  \leavevmode
  \begingroup
  \bfseries
  #2.\quad
  \endgroup
  \addcontentsline{toc}{subsection}{%
    #1}%
  }
\makeatother

\makeatletter
\def\subsubsection{\@ifstar\unnumberedsubsubsection\numberedsubsubsection}
\def\numberedsubsubsection{\@ifnextchar[
  \numberedsubsubsectionwithtwoarguments\numberedsubsubsectionwithoneargument}
\def\unnumberedsubsubsection{\@ifnextchar[
  \unnumberedsubsubsectionwithtwoarguments\unnumberedsubsubsectionwithoneargument}
\def\numberedsubsubsectionwithoneargument#1{\numberedsubsubsectionwithtwoarguments[#1]{#1}}
\def\unnumberedsubsubsectionwithoneargument#1{\unnumberedsubsubsectionwithtwoarguments[#1]{#1}}
\def\numberedsubsubsectionwithtwoarguments[#1]#2{%
  \ifhmode\par\fi
  \removelastskip
  \vskip 3ex\goodbreak
  \refstepcounter{subsubsection}%
  \noindent
  \leavevmode
  \begingroup
  \bfseries
  \S \thesubsubsection\ 
  #2.\quad
  \endgroup
  \addcontentsline{toc}{subsubsection}{%
    \protect\numberline{\thesubsubsection}%
    #1}%
  }
\def\unnumberedsubsubsectionwithtwoarguments[#1]#2{%
  \ifhmode\par\fi
  \removelastskip
  \vskip 3ex\goodbreak
  \noindent
  \leavevmode
  \begingroup
  \bfseries
  #2.\quad
  \endgroup
  \addcontentsline{toc}{subsubsection}{%
    #1}%
  }
\makeatother

\begin{document}
\begin{center}
{\Large \scaps{Extremal results regarding $K_6$-minors\\ in graphs of girth at least $5$}}\\ \ \\

{\small
Elad Aigner-Horev and Roi Krakovski\\ \ \\  

$\{$horevel,roikr$\}$@cs.bgu.ac.il\\
Department of Computer Science\\
Ben-Gurion University of the Negev,\\ 
Beer Sheva, 84105, Israel}
\end{center}

\noindent
{\small
\bfm{Abstract.} We prove that every $6$-connected graph of girth $\geq 6$ has a $K_6$-minor and thus settle Jorgensen's conjecture for graphs of girth $ \geq 6$.  
Relaxing the assumption on the girth, we prove that every $6$-connected $n$-vertex graph of size $\geq 3\frac{1}{5}n-8$ and of girth $\geq 5$ contains a $K_6$-minor.}\\

\noindent
{\small 
\scaps{Preamble.} Whenever possible notation and terminology are that of~\cite{diestel}.
Throughout, a graph is always simple, undirected, and finite. $G$ always denotes a graph.
We write $\delta(G)$ and $d_G(v)$ to denote the minimum degree of $G$ and the degree of a vertex $v \in V(G)$, respectively. $\kappa(G)$ denotes the vertex connectivity of $G$. The \emph{girth} of $G$ is the length of a shortest circuit in $G$. Finally, the cardinality $|E(G)|$ is called the \emph{size} of $G$ and is denoted $\size{G}$; $|V(G)|$ is called the \emph{order} of $G$ and is denoted $|G|$.}

\section{Introduction} A conjecture of Jorgensen postulates that {\sl the $6$-connected graphs not containing $K_6$ as a minor are the apex graphs}, where a graph is apex if it contains a vertex removal of which results in a planar graph. The $6$-connected apex graphs contain triangles. Consequently, if Jorgensen's conjecture is true, then a $6$-connected graph of girth $\geq 4$  contains a $K_6$-minor. Noting that the extremal function for $K_6$-minors is at most $4n-10$~\cite{mader2} (where $n$ is the order of the graph), our first result in this spirit is that 

\begin{statement}
\label{girth6-main}
a graph of size $\geq 3n-7$ and girth at least $6$ contains a $K_6$-minor. 
\end{statement}

So that,

\begin{statement}
every $6$-connected graph of girth $\geq 6$ contains a $K_6$-minor;
\end{statement}

\noindent
This settles Jorgensen's conjecture for graphs of girth $\geq 6$. Relaxing the assumption on the girth in \sref{girth6-main}, we prove the following.
 
\begin{statement}
\label{girth5-main}
A $6$-connected graph of size $\geq 3\frac{1}{5}n - 8$ and girth at least $5$ contains a $K_6$-minor.
\end{statement}

\noindent
\scaps{Remark.} In our proofs of \sref{girth6-main} and \sref{girth5-main}, the proofs of claims 
(\sref{girth6-main}.A-B) and (\sref{girth5-main}.A-D) follow the approach of~\cite{mader1}.

\section{Preliminaries} 
Let $H$ be a subgraph of $G$, denoted $H \subseteq G$.
The boundary of $H$, denoted by $\bnd{H}{G}$ (or simply $\bnd{H}{}$), is the set of vertices of $H$ incident with $E(G)\sm E(H)$. By $int_G H$ (or simply $int H$) we denote the subgraph induced by $V(H) \sm bnd H$. If $v \in V(G)$, then $N_H(v)$ denotes $N_G(v) \cap V(H)$.

Let $k\geq 1$ be an integer. By $k$-\emph{hammock} of $G$ we mean a connected subgraph $H\subseteq G$ satisfying $|bnd H| = k$. A hammock $H$ coinciding with its boundary is called {\sl trivial}, \emph{degenerate} if $|H| = |bnd H| +1$, and \emph{fat} if $|H|\geq |bnd H|+2$.
A proper subgraph of $H$ that is a $k$-hammock is called a \emph{proper} $k$-hammock of $H$. 
A fat $k$-hammock is called {\sl minimal} if all its proper $k$-hammocks, if any, are trivial or degenerate. Clearly,
\begin{equation}\label{each}
\mbox{every fat $k$-hammock contains a minimal fat $k$-hammock.}
\end{equation}
%

Let $H$ be a fat $2$-hammock with $bnd H = \{u,v\}$. By \emph{capping} $H$ we mean $H+uv$ if $uv \notin E(H)$ and $H$ if $uv \in E(H)$. In the former case, $uv$ is called a \emph{virtual} edge of the capping of $H$. The set $bnd H$ is called the \emph{window} of the capping.

Let now $\kappa(G)=2$ and $\delta(G) \geq 3$. By the standard decomposition of $2$-connected graphs into their $3$-connected components~\cite[Section 9.4]{bm}, such a graph has at least two minimal fat $2$-hammocks whose interiors are disjoint and that capping of each is $3$-connected. Such a capping is called an \emph{extreme} $3$-connected component. 

A $k$-(vertex)-disconnector, $k\geq 1$, is called \emph{trivial} if removal of which isolates a vertex. Otherwise, it is called \emph{nontrivial}. A graph is called \emph{essentially} $k$-connected if all its $(k-1)$-disconnectors are trivial. If each $(k-1)$-disconnector $D$ isolates a vertex and $G-D$ consists of precisely $2$ components (one of which is a singleton) then $G$ is called \emph{internally} $k$-connected.

Suppose $\kappa(G) \geq 1$ and that $D \subseteq V(G)$ is a $\kappa(G)$-disconnector of $G$.
Then, $G[C \cup D]$ is a fat $\kappa(G)$-hammock for every non-singleton component $C$ of $G-D$. 
In particular, we have that

\begin{statement}\label{fat1}
if $\kappa(G) \geq 1$, $\delta(G)\geq 3$, and $D \subseteq V(G)$ is a nontrivial $\kappa(G)$-disconnector of $G$, then $G$ has at least two fat minimal $\kappa(G)$-hammocks whose interiors are disjoint. 
\end{statement}
%

\begin{statement}\label{fat3}
If $\kappa(G) \geq 1$, $\delta(G) \geq 3$, $e \in E(G)$, and $G$ has a nontrivial $\kappa(G)$-disconnector, then $G$ has a minimal fat $\kappa(G)$-hammock $H$ such that if $e \in E(H)$, then $e$ is spanned by $bnd H$. 
\end{statement}

Let $H$ be a $k$-hammock. By \emph{augmentation} of $H$ we mean the graph obtained from $H$ by adding a new vertex and linking it with edges to each vertex in $bnd H$. 

\begin{statement}\label{fat2}
Suppose $\kappa(G)=3$ and that $H$ is a minimal fat $3$-hammock of $G$. 
Then, an augmentation of $H$ is $3$-connected.
\end{statement}

\proof
Let $H'$ denote the augmentation and let $\{x\} = V(H') \sm V(H)$. 
Assume, to the contrary, that $H'$ has a minimum disconnector $D$, $|D| \leq 2$. 
If $H'-D$ has a component containing $x$, then $H$ has a nontrivial $|D|$-hammock; contradicting the assumption that $\kappa(G)=3$. Hence, $x \in D$. 
As $x$ is $3$-valent, $H'-D$ has a component $C$ containing a single member of $bnd H'$ (= $N_{H'}(x)$), say $u$. Since $\delta(G) \geq 3$, $|N_C(u) \sm D| \geq 1$ so that $(D \sm \{x\}) \cup \{u\}$ is a disconnector of $H$ of size $\leq 2$ not containing $x$ and hence also a disconnector of $G$; contradiction. 
\QED

\begin{statement}\label{fat2'}
Suppose $\kappa(G)=3$ and that $H$ is a triangle free minimal fat $3$-hammock of $G$ such that $e \in E(G[bnd H])$. Then, an augmentation of $H-e$ is $3$-connected.
\end{statement}

\proof
Let $H'$ be the augmentation of $H-e$, let $\{x\} = V(H') \sm V(H)$, and let $e = tw$ such that $t,w \in N_{H'}(x)$. By \sref{fat2}, $\kappa(H'+e)\geq 3$. Suppose that $\kappa(H') <3$, then $H'$ contains a $2$-disconnector, say $\{u,v\}$, so that $H'=H_1 \cup H_2$, $H'[\{u,v\}] = H_1 \cap H_2$ and such that $x \in V(H_i)$ for some $i \in \{1,2\}$. Unless $x \in \{u,v\}$, then $t,w \in V(H_i)$. Thus, if $x \notin \{u,v\}$, then $\{u,v\}$ is a $2$-disconnector of $H'+e$; contradiction. 

Suppose then that, without loss of generality, $x =u$. Thus, since $x$ is $3$-valent, there exists an $i \in \{1,2\}$ such that $|N_{H_i}(x)\sm \{v\}|=1$. As $\{x,v\}$ is a minimum disconnector of $H'$, it follows that $H_i - \{x,v\}$ is connected so that $N_{H_i}(x) \cup \{
v\}$ is the boundary of a $2$-hammock of $G$; such must be trivial as $\kappa(G)=3$, implying that
$|V(H_i)| = \{x,v,z\}$, where $z \in \{t,w\}$. 

We may assume that $x$ is not adjacent to $v$; for otherwise, $|N_{H_{3-i}}(x)\sm \{v\}| =1$ so that the minimality of the disconnector $\{x,v\}$ implies that $H_{3-i}-\{x,v\}$ is connected and consequently that $N_{H_{3-i}}(x) \cup \{v\}$ is the boundary of a $2$-hammock of $G$; since such must be trivial we have that $H$ is a triangle (consisting of $\{t,v,w\}$) contradicting the assumption that $H$ is triangle-free.

Hence, since $H$ is triangle free and since each member of $\{v\} \cup N_{H_{3-i}}(x)$ 
has at least two neighbors in $H_{3-i}$, $\{v\} \cup N_{H_{3-i}}(x)$
is the boundary of a proper fat $3$-hammock of $H$; contradiction to $H$ being minimal.\QED \\

The maximal $2$-connected components of a connected graph are called its \emph{blocks}.
Such define a tree structure for $G$ whose leaves are blocks and are called the \emph{leaf} blocks of $G$~\cite{diestel}. 

We conclude this section with the following notation. Let $H \subseteq G$ be connected (possibly $H$ is a single edge). By $G/H$ we mean the contraction minor of $G$ obtained by contracting $H$ into a single vertex. We always assume that after the contractions the graph is kept simple; i.e., any multiple edges resulting from a contraction are removed.

\section{Truncations} 
Let $\F$ be a family of graphs (possibly infinite). A graph is $\F$-\emph{free} if it contains no member of $\F$ as a subgraph. A graph $G$ is \emph{nearly} $\F$-free if it is either $\F$-free or has a \emph{breaker} $x \in V(G) \cup E(G)$ such that $G-x$ is $\F$-free. A breaker that is a vertex is called a \emph{vertex-breaker} and an \emph{edge-breaker} if it is an edge.

An $\F$-\emph{truncation} of an $\F$-free graph $G$ is a minor $H$ of $G$ that is nearly $\F$-free such that either $H \subseteq G$ (and then it has no breaker) or $H$ contains a breaker $x$ such that $H-x \subseteq G$. In the former case, the truncation is called \emph{proper}; in the latter case, the truncation is \emph{improper} with $x$ as its breaker and $H-x$ as its \emph{body}. 
An improper truncation is called an \emph{edge-truncation} if its breaker is an edge and a \emph{vertex-truncation} if its breaker is a vertex. A vertex-truncation is called a \emph{$3$-truncation} if its breaker is $3$-valent. 
 
\begin{statement}
\label{ess4}
Let $\F$ be a graph family such that $K_3 \in \F$ and let $G$ be $\F$-free with $\delta(G) \geq 3$. 
Then $G$ has an essentially $4$-connected $\F$-truncation $H$ such that: 

(\sref{ess4}.1) $|H| \geq 4$; and

(\sref{ess4}.2) if $H$ is a vertex-truncation then it is a $3$-truncation and $|H| \geq 5$.
\end{statement}

\proof Let $\H$ denote the $3$-connected truncations of $G$.\\

\noindent
(\sref{ess4}.A) \emph{$\H$ is nonempty. In particular, $\H$ contains a truncation $H$ with $|H| \geq 4$ so that if improper then it is an edge-truncation with edge-breaker $e$ such that $\kappa(H-e)=2$.}\\

\noindent
\emph{Subproof.}
We may assume that $G$ is connected. Let $B$ be a leaf block of $G$ (possibly $B=G$). If $\kappa(B) \geq 3$, then (\sref{ess4}.1) follows (by setting $H=B$) as $B$ is a proper truncation of $G$. Assume then that $\kappa(B) = 2$ and let $H$ be an extreme $3$-connected component of $B$ with window $\{x,y\}$. Now, $H \in \H$ with possibly $xy$  an edge-breaker. If $H$ is improper, then $\kappa(H-xy) =2$. Note that $\delta(G) \geq 3$ implies that $|H| \geq 4$ in both cases. 
\inQED \\

If $\H$ contains a proper or an edge-truncation that is essentially $4$-connected, then (\sref{ess4}.1) follows. Suppose then that 
\begin{equation}\label{assume}
\mbox{$\H$ has no proper or edge-truncations that are essentially $4$-connected.}
\end{equation}

\noindent
(\sref{ess4}.B) \emph{Assuming (\ref{assume}), then $\H$ contains a truncation that if improper then it is a $3$-truncation of order $\geq 5$.}\\

\noindent
\emph{Subproof.} Let $H \in \H$ such that if improper then $H$ and $e$ are as in (\sref{ess4}.A). By (\ref{assume}) and \sref{fat3}, $H$ has a minimal fat $3$-hammock $H'$ such that if $e \in E(H')$, then $e$ is spanned by the boundary of $H'$. Let $H''$ be the graph obtained from an augmentation of $H'$ by removing $e$ if it is spanned by $bnd H'$. Let $\{x\} = V(H'') \sm V(H')$.

By \sref{fat2} and \sref{fat2'}, $\kappa(H'') \geq 3$ so that $H'' \in \H$ with $x$ as a potential $3$-valent vertex-breaker and (\sref{ess4}.B) follows. 

%
Finally, note that $|int H'| \geq 2$ so that $|H''| \geq 5$. 
\inQED \\

Next, we show the following.\\

\noindent
(\sref{ess4}.C) \emph{If $\H$ contains a $3$-truncation $X$ of order $\geq 5$, then $\H$ contains essentially $4$-connected $3$-truncations $Y$ such that $5 \leq |Y| \leq |X|$.}\\

\noindent
\emph{Subproof.}
Let $H^* \in \H$ be a $3$-truncation of order $\geq 5$ with the order of its body minimized. We show that $H^*$ is essentially $4$-connected. Let $x$ denote the vertex-breaker of $H^*$. By the minimality of $H^*$, 
\begin{equation}
\label{red}
\mbox{any minimal fat $3$-hammock $T$ of $H^*$ with $x \notin V(T)$ satisfies $T = H^*-x$}
\end{equation}
(so that $bnd T = N_{H^*}(x)$). 

Assume now, towards contradiction, that $H^*$ is not essentially $4$-connected so that it contains nontrivial $3$-disconnectors and at least two minimal fat $3$-hammocks that may meet only at their boundary, by \sref{fat1}. By (\ref{red}), existence of at least two such hammocks implies that $x$ belongs to every nontrivial $3$-disconnector and thus to the boundary of every minimal fat $3$-hammock. As $x$ is $3$-valent, there is a minimal fat $3$-hammock $T$ of $H^*$ with $x$ on its boundary such that $N_{T}(x)=\{y\}$. As $T$ is a minimal fat $3$-hammock, $V(T)$ consists of $x,y$, the two members of $bnd T \sm \{x\}$, and an additional vertex $u$. As $\delta(G) \geq 3$, $uy \in E(T)$, $u$ is adjacent to both members of $bnd T \sm \{x\}$ and $y$ is adjacent to at least one member of $bnd T \sm \{x\}$. Hence, $K_3 \subseteq T-x \subseteq H^*-x$ so that $x$ is not a breaker; contradiction.
\inQED\\

Assuming (\ref{assume}), then, by (\sref{ess4}.B), there are $3$-connected $3$-truncations of $G$ of order $\geq 5$ so that an essentially $4$-connected $3$-truncation of $G$ exists by (\sref{ess4}.C). \QED \\ 

\begin{statement}\label{ess4'}
Let $\F$ be a graph family such that $\{K_3,K_{2,3}\} \subseteq \F$, then $G$ has an internally $4$-connected $\F$-truncation satisfying (\sref{ess4}.1-2) and if such is a vertex-truncation then it is a $3$-truncation. 
\end{statement}

\proof Let $\T$ denote the essentially $4$-connected truncations of $G$ that are either proper, or edge-truncations, or $3$-truncations; $\T$ is nonempty by \sref{ess4}. Let $\alpha(\T)$ denote the least $k$ such that $\T$ contains a proper truncation of order $k$ or an improper edge-truncation of order $k$. Let $\beta(\T)$ denote the least $k$ such that $\T$ contains an improper $3$-truncation with its body of order $k$. Let $H \in \T$ such that $|H| = \min\{\alpha(\T),\beta(\T)+1\}$ and let $x$ denote its breaker if improper. 

We show that $H$ is internally $4$-connected. To see this, assume, to the contrary, that $H$ is not internally $4$-connected and let $D$ be a $3$-disconnector of $H$ such that $H-D$ consists of $\geq 3$ components at least one of which is a singleton (since $H$ is essentially $4$-connected).
Let $\C$ denote the non-singleton components of $H-D$. Since $K_{2,3} \in \F$, $|\C| \geq 1$ 
 
Suppose $J = H[C \cup D]$ is a $3$-hammock of $H$, for some $C \in \C$, that does not meet $x$ in its interior (if $x$ exists). By the choice of $H$, 
\begin{equation}\label{on}
\mbox{for each fat $3$-hammock $X$ of $J$ either $x \in bnd X$ or $x \in E(H[bnd X])$.}
\end{equation}
Indeed, for otherwise, an augmentation of a minimal fat $3$-hammock of $X$ is a $3$-truncation of order $\geq 5$ of $G$ that belongs to $\H$ and has order $< |H|$, where $\H$ is as in the proof of \sref{ess4}; existence of such a $3$-truncation of $G$ implies that $G$ has an essentially $4$-connected $3$-truncation of order $\geq 5$, by (\sref{ess4}.C), and such has order $< |H|$  contradicting the choice of $H$. Consequently, the assumption that the interior of $J$ does not meet $x$ implies that 
\begin{equation}\label{on2}
\mbox{if $J$ exists, then $x\in D \cup E[H[D]]$.} 
\end{equation}

Suppose now that $J$ has a minimal fat $3$-hammock $J'$ (possibly $J'=J$) with $x \in bnd J'$ so that $x \in D$, by (\ref{on2}). $|D| = \kappa(H)$ imply that $x$ is incident with each component of $H-D$ so that $|N_{int J'}(x)|=1$, as $x$ is $3$-valent. The minimality of $J'$ then implies that $|int J'|=2$ so that $J'-x$ contains a $K_3$ (see proof of (\sref{ess4}.C) for the argument) and thus $x$ is not a breaker of $H$; contradiction.
 
Suppose next that $J'$ is a minimal fat $3$-hammock of $J$ whose boundary vertices span $x$ (as an edge). Then, an augmentation of $J'-x$ belongs to $\H$, by \sref{fat2'}, and such contains an essentially $4$-connected $3$-truncation of $G$, by (\sref{ess4}.C), of order $< |H|$. Hence,
\begin{equation}\label{on4}
\mbox{$J$ (if exists) has no minimal fat $3$-hammock $J'$ with $x \in bnd J' \cup E[H[bnd J']]$.}
\end{equation}

If $J$ exists, then (\ref{on}) and (\ref{on4}) are contradictory. Thus, to obtains a contradiction and hence conclude the proof of \sref{ess4'} we show that a $3$-hammock such as $J$ exists. This is clear if $|\C| \geq 2$ as then at least one member of $\C$ does not meet $x$. Suppose then that $|\C|=1$ so that $H-D$ consists of two singleton components, say $\{u,v\}$, and the single member $C$ of $\C$. $D \cup \{u,v\}$ induce a $K_{2,3}$, say $K$. Since $K_{2,3} \in \F$ and $x$ is a breaker, $K$ contains $x$ so that $C$ does not; hence, $H[C \cup D]$ is the required $3$-hammock.\QED \\

For $k \geq 4$, a graph that is nearly $\{K_3,C_4,\ldots,C_{k-1}\}$-free is called \emph{nearly $k$-long}. That is, $G$ is nearly $k$-long if either it has girth $\geq k$ or it has a \emph{breaker} $x \in V(G) \cup E(G)$ such that $G-x$ has girth $\geq k$. 

A nearly $5$-long graph is nearly $\{K_3,C_4\}$-free; such is also nearly $\{K_3,K_{2,3}\}$-free. 
In addition, a $3$-connected nearly $5$-long truncation has order $\geq 5$. Consequently, we have the following consequence of \sref{ess4'}.

\begin{statement}\label{trun}
A graph with girth $\geq k \geq 5$ and $\delta \geq 3$ has an internally $4$-connected nearly $k$-long truncation of order $\geq 5$ and if such is a vertex-truncation then it is a $3$-truncation. 
\end{statement}

\section{Nearly long planar graphs} 
For a plane graph $G$, we denote its set of faces by $F(G)$ and by $X_G$ its infinite face. 

\begin{statement}
\label{dis:girth6}
Let $G$ be a $2$-connected plane graph of girth $\geq 6$, and let $S \subseteq V(G)$ be the $2$-valent vertices of $G$. Then, $|S| \geq 6$. 
\end{statement}

\proof By Euler's formula:
\begin{equation}
\label{eq:1}
\mbox{ $|E(G)|=|V(G)|+ |F(G)|-2$.}
\end{equation} 
Since $G$ is $2$-connected, every vertex in $V(G)\setminus S$ is at least $3$-valent so that
\begin{equation}
\label{eq:2}
\mbox{ $2|E(G)|\geq 3(|V(G)|-|S|)+ 2|S|$.}
\end{equation} 
As $G$ is of girth $\geq 6$ and $2$-connected (and hence every edge is contained in exactly two distinct faces) then:
\begin{equation}
\label{eq:3}
\mbox{ $2|E(G)|\geq 6|F(G)|$.}
\end{equation} 
Substituting (\ref{eq:1}) in (\ref{eq:2}),
\begin{equation}
\label{eq:4}
\mbox{ $2(|V(G)|+ |F(G)|-2)\geq 3(|V(G)|-|S|)+ 2|S|  \Rightarrow  |V(G)| \leq 2|F(G)|+|S|-4$}
\end{equation}
Substituting (\ref{eq:1}) in (\ref{eq:3}),
\begin{equation}
\label{eq:5}
\mbox{ $2(|V(G)|+ |F(G)|-2)\geq 6|F(G)|  \Rightarrow  |V(G)|\geq 2|F(G)|+2$}
\end{equation}

From (\ref{eq:4}) and (~\ref{eq:5}),
 \begin{equation}
\label{eq:6}
\mbox{ $2|F(G)|+2 \leq 2|F(G)|+|S|-4 \Rightarrow |S| \geq 6$}
\end{equation}
Hence, the proof follows. \QED \\

From~\sref{dis:girth6} we have that:

\begin{statement}
\label{dis2}
A nearly $6$-long internally $4$-connected graph is nonplanar. 
\end{statement}

\begin{statement}
\label{sizeofnearlylong}
Let $G$ be a nearly $5$-long internally $4$-connected planar graph and suppose that if $G$ has a vertex-breaker, then it also has a vertex-breaker which is a $3$-valent vertex. Then, $|G|\geq 11$. 
\end{statement}
\proof Define $S \subseteq V(G)\cup E(G)$ as follows. If $G$ is of girth $\geq 5$ set $S:=\emptyset$; otherwise set $S:=\{x\}$, where $x \in V(G)\cup E(G)$ is a breaker of $G$ so that if $x \in V(G)$ then $x$ is $3$-valent. Then, $G- S$ is $2$-connected, and has at most three $2$-valent vertices. Hence,
\begin{equation}
\label{eq:12}
\mbox{ $2|E(G)|\geq 3(|V(G)|-3)+ 6$.}
\end{equation} 
As $G- S$ is of girth $\geq 5$ and $G$ is $2$-connected then:
\begin{equation}
\label{eq:13}
\mbox{ $2|E(G)|\geq 5|F(G)|$.}
\end{equation} 
Substituting (\ref{eq:1}) in (\ref{eq:12}),
\begin{equation}
\label{eq:14}
\mbox{ $2(|V(G)|+ |F(G)|-2)\geq 3(|V(G)|-3)+ 6  \Rightarrow  |F(G)| \leq (|V(G)|+1)/2$}
\end{equation}
Substituting (\ref{eq:1}) in (\ref{eq:13}),
\begin{equation}
\label{eq:15}
\mbox{ $2(|V(G)|+ |F(G)|-2)\geq 5|F(G)|  \Rightarrow  |F(G)|\geq (2|V(G)|-2)/3$}
\end{equation}
From (\ref{eq:14}) and (\ref{eq:15}),
\begin{equation}
\label{eq:16}
\mbox{ $(|V(G)|+1)/2 \leq (2|V(G)|-2)/3 \Rightarrow |V(G)|\geq 11$}
\end{equation}
Hence, the proof follows. \QED

\begin{statement}
\label{dis1}
A $2$-connected plane graphs $G$ satisfying the following does not exist.

(\sref{dis1}.1) $G$ has girth $\geq 5$;  

(\sref{dis1}.2) each member of $V(G)- V(X_G)$ is at least $4$-valent; and
	
(\sref{dis1}.3) $G$ has a set $S \subseteq V(X_G)$, $|S| \leq 3$ (possibly $S= \emptyset$) with each of its members \indent \indent \indent \ \ $2$-valent and each member of $V(X_G)- S$ at least $3$-valent. 

\end{statement}

\proof
Assume towards contraction that the claim is false. We will use the Discharging Method to obtain a contradiction to Euler's formula.  The discharging method starts by assigning numerical values (known as charges) to the elements of the graph. For $x \in V(H) \cup F(H)$, define $ch(x)$ as follows.\\

(CH.1) $ch(v)=6-d_H(v)$, for any $v \in V(H)$.

(CH.2)  $ch(f)=6-2|f|$, for any $f \in F(H) - \{X_{H}\}$.

(CH.3) $ch(X_{_H})=-5\frac{2}{3}-2|X_{_H}|$.\\

\noindent
Next, we show that

\begin{equation}\label{total}
\displaystyle\sum_{x \in V(H) \cup F(H)} ch(x) = \frac{1}{3}.
\end{equation}

\noindent
{\sl Proof.} 
\begin{eqnarray*}
\displaystyle\sum_{x \in V(H) \cup F(H)} ch(x) &=& -5\frac{2}{3} -2|X_H|+  \displaystyle\sum_{f \in F(H)- X_H} (6-2|f|)  + \displaystyle\sum_{v \in V(H)} (6-d(v))\\ &=& -5\frac{2}{3} -2|X_H| + 6(|f(H)|-1) + \displaystyle\sum_{f \in F(H)- X_H}(-2|f|) + \displaystyle\sum_{v \in V(H)} (6-d(v))\\
                &=& -5\frac{2}{3}+6(|f(H)|-1)-2(2|E|)+6|V(H)|-2|E(H)|\\
                &=& 6(F(H)-E(H)+V(H))-11\frac{2}{3}=\frac{1}{3}
\end{eqnarray*}

\noindent
Next the charges are locally redistributed according to the following discharging rules: 

\begin{itemize}
\item [(DIS.1)] If $v$ is 2-valent, then $v$  sends $3\frac{1}{5}$ to $X_G$ and $\frac{4}{5}$ to the other face incident to it.

\item [(DIS.2)] If $v$ is 3-valent, then $v$  sends $1\frac{5}{8}$ to $X_G$ and $\frac{4}{5}$ to every other face incident to it.

\item [(DIS.3)] If $v$ is at least 4-valent, then $v$ sends $\frac{4}{5}$ to each incident face.

\end{itemize} 

For $x \in V(G) \cup F(G)$, let $ch^*(x)$ (denoted as the modified charge) be the resultant charge
after modification of the initial charges according to (DIS.1-3). 
We obtain a contradiction to (\ref{total}) by showing that $ch^*(x) \leq 0$ for every $x \in V(H) \cup F(H)$. This is clearly implied by the following claims proved below.\\

\noindent
(A) \emph{$ch^*(v) \leq 0$, for each $v \in V(H)$.}\\
(B) \emph{$ch^*(f) \leq 0$, for each $f \in F(H) - \{X_{H}\}$.}\\
(C) \emph{$ch^*(X_{_H}) \leq 0$.}\\

Observe that according to DIS.(1)-(3), faces do not send charge and vertices do not receive charge.\\ 

\noindent
{\sl Proof of} (A). It is sufficient to consider vertices $v$ satisfying 
$d_G(v) \geq 5$. Indeed, if $d_H(v) \geq 6$, then $ch(v)=ch^*(v) \leq 0$ by (CH.1). 
If $2 \leq d_G(v) \leq 3$, then it is easily seen by (CH.1) and (DIS.1-2) that $ch^*(v)=0$.
If $4 \leq d_G(v) \leq 5$, then, by (CH.1) and (DIS.3), $ch^*(v)=6-d_H(v)- \frac{4}{5}d_G(v) \leq 0$. \inQED\\
 
\noindent
{\sl Proof of} (B). Let $f \in F(H) - \{X_{_H}\}$. By (DIS.1-3), $f$ receives a charge of $\frac{4}{5}$ from every vertex incident to it. Hence, togther with (CH.2), $ch^*(f)=6-2|f| + \frac{4}{5}|f| \leq 0$. (The last inequality follows as $|f| \geq 5$.) \inQED. \\

\noindent
{\sl Proof of} (C). Let $S_1 \subseteq V(X_G)$ be the set of 3-valent vertices of $X_G$, and let $S_2=V(X_G) - (S \cup S_1)$. By (CH.3), (DIS.1-3) and as $|S| \leq 3$, we see that $ch^*(f)=-5\frac{2}{3}-2|X_{G}| +3\frac{1}{5}|S|+1\frac{5}{8}|S_1|+ \frac{4}{5}|S_2| \leq -5\frac{2}{3}-2|X_{G}| + 3 \times 3\frac{1}{5} +1\frac{5}{8}(|X_G|-3) = -\frac{3}{8}|X_G|-\frac{11}{12} \leq 0$. \inQED
\QED

\section{$K_5$-minors in internally 4-connected graphs}
By $V_8$ we mean $C_8$ together with $4$ pairwise overlapping chords. 
By $TG$ we mean a subdivided $G$.

The following is due to Wanger.

\begin{statement} \bfm{\cite[Theorem 4.6]{wagner}}
\label{transform-V8}
If $G$ is $3$-connected and $TV_8 \subseteq G$ then either $G \iso V_8$ or $G$ has a $K_5$-minor. 
\end{statement}

The following structure theorem was proved independently by Kelmans~\cite{kel} and Robertson~\cite{Robertson}.

\begin{statement}\bfm{\cite{kel}}
\label{robertson}
Let $G$ be internally $4$-connected with no minor isomorphic to $V_8$. Then
$G$ satisfies one of the following conditions:

(\sref{robertson}.1) $G$ is planar; 

(\sref{robertson}.2) $G$ is isomorphic to the line graph of $K_{3,3}$;

(\sref{robertson}.3) there exist a $uv \in E(G)$ such that $G-\{u,v\}$ is a circuit;

(\sref{robertson}.4) $|G| \leq 7$;

(\sref{robertson}.5) there is an $X \subseteq V(G)$, $|X| \leq 4$ such that $\size{G - X} =0$.
\end{statement}

From~\sref{transform-V8} and~\sref{robertson} we deduce that

\begin{statement}
\label{thomas-cor}
A nearly $5$-long internally $4$-connected nonplanar $G$ has a $K_5$-minor.
\end{statement}

\proof
We may assume that $G \notiso V_8$ and that $G$ has no $V_8$-minor. The former since $V_8$ is not nearly $5$-long and the latter by \sref{transform-V8}. Hence, $G$ satisfies one of (\sref{robertson}.1-5). As $G$ is nonplanar, by assumption, and the line graph of $K_{3,3}$ has a $K_5$-minor (and is not nearly $5$-long) it follows that $G$ satisfies one of (\sref{robertson}.3-5).

If $G$ is of girth $\leq 4$, let $a \in V(G)\cup E(G)$ be a breaker of $G$; otherwise (if $G$ has girth $\geq 5$) let $a$ be an arbitrary vertex of $G$. If $a \in V(G)$, put $b:=a$; otherwise let $b$ be some end of $a$. By defintion, $G-b$ has girth $\geq 5$.\\

\noindent
(\sref{thomas-cor}.A) \emph{$G- \{u,v\}$ is not a circuit for any $u,v \in V(G)$ so that $G$ does not satisfy (\sref{robertson}.3)}. \\

\noindent
\emph{Subproof.} For suppose not; and let $C:=G - \{u,v\}=\{x_0,\dots,x_{k-1}\}$, where $k\geq 3$ is an integer. 

Suppose first that $b\in \{u,v\}$ and assume, without loss of generality, that $u=b$. Then, $k\geq 5$. As $v$ is at least $3$-valent, there exists $0 \leq i \leq k-1$ so that $vx_i \in E(G)$. Since $G-b$ has girth $\geq 5$, $vx_{i+1},vx_{i+2} \not\in E(G)$ (subscript are read modulo $k$). Since $x_{i+1}$ and $x_{i+2}$ are at least $3$-valent in $G$, each is adjacent to $u$. But then $\{u,x_{i},x_{i+3}\}$ is a $3$-disconnector of $G$ separating $\{x_{i+1},x_{i+2}\}$ from $\{v,x_{i+4}\}$ (note that since $k\geq 5$,  $x_{i+1},x_{i+2} \neq x_{i+4}$); a contradiction to $G$ being internally $4$-connected.

Suppose then that $x_i=b$, for some $0 \leq i \leq k-1$. Hence, exactly one of $v$ and $u$ is adjacent to $x_{i+1}$ and exactly one to $x_{i+2}$ (this is true since every vertex of $C$ is adajcent to $v$ or $u$, and if say, $v$, is adajcent to both $x_{i+1}$ and $x_{i+2}$ then $G-b$ conatins a trinagle). If $x_{i+3} \neq x_i$, then $x_{i+3}$ is adjacent to one of $u$ and $v$. If $x_i=x_{i+3}$, then $C$ is a circuit of length three, and $V(G)=5$. Both cases contradict the fact that $G$ is nearly $5$-long.\inQED \\

\noindent
(\sref{thomas-cor}.B) \emph{$|G|\geq 8$ so that $G$ does not satisfy (\sref{robertson}.4). } \\

\noindent
\emph{Subproof.} For suppose $|G| \leq 7$. As $G$ is internally $4$-connected, $G- b$ is $2$-connected. Since $G-b$ is of girth $\geq 5$, then $G-b$ contains an induced circuit $C$ of length $\geq 5$. Hence $|G| \geq 6$. If $|G|=6$, then $G=C \cup b$ and then $G$ is planar; a contracation. If $|G|=7$ then $G$ is a circuit plus two vertices and we get a contrdaction to (\sref{thomas-cor}.A).  Hence, $V(G) \geq 8$.\inQED \\

To reach a contradiction we show that (\sref{robertson}.5) is not satisfied by $G$. For suppose it is satisfied and let $X$ be as in (\sref{robertson}.5) and let $Y=V(G)- X$. As $V(G) \geq 8$, then $|Y|\geq 4$ and every vertex of $Y$ is adjacent to at least three vertices in $X$. But then it is easily  seen that $G$ is of girth $\leq4$ but contains no edge- or vertex-breaker; a contradiction. \QED \\

Let $G$ be a plane graph. By \emph{jump} over $G$ we mean a path $P$ internally-disjoint of $G$ whose ends are not cofacial in $G$. 

\begin{statement}
\label{jump}
Let $G$ be an internally $4$-connected nearly $5$-long plane graph and let $P$ be a jump over $G$. Then, $G$ has a $K_5$-minor with every branch set meeting $V(G)$.
\end{statement}

\proof Put $G':= G \cup P$. (By possibly contracting $P$) we may assume that $P$ is an edge $e$ with both ends in $G$.  Suffices now to show that $G'$ has a $K_5$-minor. Suppose $G'$ has no such minor. We may assume that $G' \notiso V_8$, since $V_8$ with any edge removed is not internally $4$-connected, and that $G'$ has no $V_8$-minor, by \sref{transform-V8}. Since $G'$ is nonplanar, $|G'| \geq |G|\geq 11$, by~\sref{sizeofnearlylong}, and since the line graph of $K_{3,3}$ has a $K_5$-minor, we have that $G'$ satisfies (\sref{robertson}.3) or (\sref{robertson}.5). We show that both options lead to a contradiction to the definition of $G$.  

Suppose (\sref{robertson}.3) is satisfied.  Set $C:=G'- \{u,v\}=\{x_0,\dots,x_{k-1}\}$, where $k\geq 9$ is an integer. If $e \not\in E(C)$, then a contradiction is obatined by showing that $G-e-\{v,u\}$ cannot be  a circuit. The proof is exactly the same as the proof of (\sref{thomas-cor}.A) with $G-e$ instead of $G$.

Hence we may assume that $e \in E(C)$; so let $e=x_ix_{i+1}$, for some $0\leq i \leq k-1$ (subscript are read modulo $k$). Observe that $d_{G'}(x_i),d_{G'}(x_{i+1}) \geq 4$. Hence, in $G$, each of $x_i$ and $x_{i+1}$ is adajcent to both $u$ and $v$. 

By assumtion that (\sref{robertson}.3) is satisfied, $uv \in E(G)$, and we see that one of $u$ or $v$ is a breaker, say $u$. Hence, $vx_{i+2},vx_{i+3} \notin E(G)$. But then, since and $d_G(x_{i+1}),d_G(x_{i+2})= 3$, the set $\{u,x_{i+1},x_{i+4}\}$  is a 3-disconnector of $G$ (note that since $k\geq 9$, $x_{i+1},x_{i+4}$ are distinct) separating $\{x_{i+2},x_{i+3}\}$ from $\{x_{i+5},x_{i+6}\}$; a contradiction. Hence (\sref{robertson}.3) is not satisfied.

Suppose (\sref{robertson}.5) is satisfied. As $V(G)\geq 11$, it is easily seen that $G ~(=G'- e)$ is of girth $\leq 4$ but has no edge- or vertex-breaker; a contradiction. This concludes the proof. \QED\\

By \emph{society} we mean a pair $(G,\Omega)$ consisting of a graph $G$ and a cyclic permutation $\Omega$ over a finite set $\overline{\Omega} \subseteq V(G)$. Let $\overline{\Omega}=\{v_1,\ldots,v_k\}$, $k \geq 4$. 
Two pairs of vertices $\{s_1,t_1\} \subseteq \overline{\Omega}$ and $\{s_2,t_2\} \subseteq \overline{\Omega}$ are said to \emph{overlap} along $(G, \Omega)$ if $\{s_1,s_2,t_1,t_2\}$ occur in $\overline{\Omega}$ in this order along $\Omega$.  

Two vertex disjoint paths $P$ and $P'$ of $G$ that are both internally-disjoint of $\overline{\Omega}$ are said to form a \emph{cross} on $(G,\Omega)$ if their ends are in $\overline{\Omega}$ and these overlap along $(G,\Omega)$. 

\begin{statement}\bfm{\emph{\cite[Lemma (2.4)]{RST}}}
\label{cross}
Let $(G,\Omega)$ be a society. Then either

(\sref{cross}.1) $(G,\Omega)$ admits a cross in $G$, or

(\sref{cross}.2) $G = G_1 \cup G_2$, $G_1 \cap G_2 = G[D]$, $|D| \leq 3$ such that
$\overline{\Omega} \subseteq V(G_1)$ and\\ \indent \indent \indent $|V(G_2) \setminus V(G_1)| \geq 2$, or

(\sref{cross}.3) $G$ can be drawn in a disc with $\overline{\Omega}$ on the boundary in order $\Omega$.
\end{statement}

Let $C$ be a circuit in a plane graph $G$. Then the clockwise ordering of $V(C)$ induced by the embedding of $G$ defines a cyclic permutation on $V(C)$ denoted $\Omega_C$ and we do not distinguish between the cyclic shifts of this order. Then, $(G,\Omega_C)$ is a society with $\overline{\Omega_C} = V(C)$. Throughout, we omit this notation when dealing with such societies of circuits of plane graphs and instead say that $C$ is a society of $G$. 

\begin{statement}
\label{addcross}
Let $G$ be a $3$-connected plane graph of order $\geq 5$ and let $P$ and $P'$ be vertex disjoint paths that are internally-disjoint of $G$ and whose ends are contained in a facial circuit $f$ of $G$. If $P\cup P'$ form a cross on $f$, then $G\cup P \cup P'$ contains a $K_5$-minor with every branch set meeting $V(G)$. 
\end{statement}

\proof
Clearly, $V(G) \not= V(f)$. Since the facial circuits of a $3$-connected plane graph are it induced nonseparating circuits~\cite{mt}, we have that $G-V(f)$ is connected so that $f \cup P \cup P'$ have a $K_4$-minor which is completed into a $K_5$-minor by adding a fifth branch set that is $G-V(f)$ (as $f$ is an induced circuit). 
\QED

\section{Proof of \sref{girth6-main}} Let $\H = \{H \subseteq G: \mbox{$H$ is connected, $|G/H| \geq 5$, and $\size{G/H}\geq 3|G/H|-7$}\}$.
$\H$ contains every member of $V(G)$ as a singleton and thus nonempty. Let $H_0 \in \H$ be maximal
in $(\H,\subseteq)$, $H_1 = G[N_G(H_0)]$, and let $G_0 = G/H_0$, where $z_0 \in V(G_0)$ represents $H_0$. Let $G_1 = G_0 -z_0$ and note that $G_1 \subseteq G$. 

$|G_0| = 5$ implies that $\size{G_0} \geq 8$ so that $\size{G_1}\geq 4$ and contains a $k$-circuit with $k< 5$; contradiction to the assumption that $G$ has girth at least $6$. Thus, we may assume that\\ 

\noindent
(\sref{girth6-main}.A) $|G_0| \geq 6$.\\ 

Let $x\in V(H_1)$ and put $G'_0 = G_0 /z_0x$. $|G'_0| \geq 5$, by (\sref{girth6-main}.A). Thus, the maximality of $H_0$ in $(\H,\subseteq)$ implies that $\size{G'_0}\leq 3|G'_0|-8$. 
Thus, $\size{G_0} - \size{G'_0} \geq 3|G_0| -7 -3(|G_0|-1)+8 \geq 4$; implying that $z_0x$ is common to at least three triangles so that $d_{H_1}(x) \geq 3$. It follows then that\\   

\noindent
(\sref{girth6-main}.B) $\delta(H_1) \geq 3$.\\ 

Let $H$ be an internally $4$-connected nearly $6$-long truncation of $H_1$, by \sref{trun}. Such is nonplanar by \sref{dis2} and has a $K_5$-minor by \sref{thomas-cor}. Consequently, $G_0$ has a $K_6$-minor. \QED

\section{Proof of \sref{girth5-main}} In a manner similar to that presented in the proof of \sref{girth6-main}, let $\H = \{H \subseteq G: \mbox{$H$ is connected, $|G/H| \geq 5$, and $\size{G/H}\geq 3\frac{1}{5}|G/H|-8$}\}$ (such is nonempty) and let $H_0,H_1,G_0,z_0,G_1$ be as in the proof of \sref{girth6-main}.

$|G_0| = 5$ implies that $\size{G_0} \geq 8$ so that $\size{G_1}\geq 4$ and contains a $k$-circuit with $k< 5$; contradiction to the assumption that $G$ has girth at least $5$. Thus, we may assume that\\ 

\noindent
(\sref{girth5-main}.A) $|G_0| \geq 6$.\\

Let $x\in V(H_1)$ and put $G'_0 = G_0 / z_0x$. $|G'_0| \geq 5$, by (\sref{girth5-main}.A). Thus, the maximality of $H_0$ in $(\H,\subseteq)$ implies that $\size{G'_0}\leq 3\frac{1}{5}|G'_0|-9$. 
Thus, $\size{G_0} - \size{G'_0} \geq 3\frac{1}{5}|G_0| -8 -3\frac{1}{5}(|G_0|-1)+9 \geq 4$; implying that $z_0x$ is common to at least three triangles so that $d_{H_1}(x) \geq 3$. It follows then that\\   

\noindent
(\sref{girth5-main}.B) $\delta(H_1) \geq 3$;\\ 

\noindent
implying that\\

\noindent
(\sref{girth5-main}.C) $\delta(G_0) \geq 4$.\\

Next, we prove that\\ 

\noindent
(\sref{girth5-main}.D) $\kappa(G_0) \geq 5$.\\ 

To see (\sref{girth5-main}.D), let $T \subseteq V(G)$ be a minimum disconnector of $G_0$ and assume, towards contradiction, that $|T| \leq 4$. As $\kappa(G) \geq 6$, $z_0 \in T$. Let then $y = |N_{G_0}(z_0) \cap T|$ and let $\C$ denote the components of $G_0-T$. Choose $C \in \C$ and put $H_1 = G_0[C \cup T]$ and $H_2 = G_0-C$. 

Let $H'_i$ be the graph obtained from $G_0$ by contracting $H_{3-i}$ into $z_0$ (note that minimality of $T$ implies that each of its members is incident with each member of $\C$), for $i=1,2$. As $|H_i| \geq 5$, by (\sref{girth5-main}.C), then $|H'_i| \geq 5$, for $i=1,2$. The maximality of $H_0$ in $(\H,\subseteq)$ then implies that $\size{H'_i}\leq 3\frac{1}{5}|H'_i|-9$.

As $z_0x \in E(H'_i)$ for each $x \in T'=T\sm\{z_0\}$, for $i=1,2$, it follows that
\begin{equation}
\size{G_0}+y+2(|T'|-y)+\size{G_0[T']} \leq \size{H'_1}+\size{H'_2} \leq 3\frac{1}{5}(|G_0|+|T|)-18.
\end{equation}
As $\size{G_0} \geq 3\frac{1}{5}|G_0| - 8$, we have that
\begin{equation}
\label{eq}
8 + \size{G_0[T']} \leq 1\frac{1}{5}|T|+y.
\end{equation}
Now, $|T| \leq 4$ (by assumption), so that $y \leq 3$, and $\size{G_0[T']} \geq 0$. Consequently, the right hand size of (\ref{eq}) does not exceed $7.8$. This contradiction establishes (\sref{girth5-main}.D).

Let $\B$ denote the bridges of $H_1$ in $G_1$. 
We may assume that $\B$ is nonempty. Otherwise, $G_1$ coincides with $H_1$ so that $H_1$ is a nonplanar $4$-connected graph of girth $\geq 5$ and thus containing a $K_5$-minor by \sref{thomas-cor}. Consequently, $G_0$ has a $K_6$-minor and \sref{girth5-main} follows.

Let $H$ be an internally $4$-connected nearly $5$-long truncation of $H_1$, by \sref{trun}. We may assume that $H$ is planar for otherwise $H$ has a $K_5$-minor, by \sref{thomas-cor}, so that $G_0$ has a $K_6$-minor and \sref{girth5-main} follows.
Let $x$ denote the breaker of $H$, if such exists in $H$. 
Let $\B_1 = \emptyset$ if $x$ does not exist (so that $H \subseteq G$) or is an edge-breaker. Otherwise (i.e., if $x$ is a vertex-breaker), $\B_1$ denotes the members of $\B$ with attachment vertices in the subgraph of $H_1$ contracted into $x$. Put $\B_2 = \B \sm \B_1$.    

Fix an embedding of $H$ in the plane. No member of $\B$ defines a jump over $H$ for otherwise the union of $H$ and such a jump has has a $K_5$-minor with every branch set meeting $V(H)$, by \sref{jump}. Hence, every member of $\B$ has all of its attachment vertices confined to a single face of $H$. 

By \emph{patch} we mean a face $f$ of $H$ together with all members of $\B$ attaching to $V(f)$. Patches not meeting $x$ in case it is a vertex-breaker are called \emph{clean} (so that if $x$ does not exist or is an edge-breaker, then every patch is clean). $f$ is called the \emph{rim} of the patch. If $\P$ is a patch with rim $f$, then by $(\P,\Omega_f)$ we mean a society with $\overline{\Omega_f} = V(f)$ and $\Omega_f$ is the clockwise order on $V(f)$ defined by the embedding of $f$ in the plane.\\

\noindent
(\sref{girth5-main}.E) \emph{Let $H'$ denote the union of $H$ and all members of $\B_2$. Then, $H'$ is planar.}\\

To see (\sref{girth5-main}.E) it is sufficient to show that every clean patch is planar. 
Indeed, since any two faces of $H$ meet either at a single vertex or at a single edge, the union of any number of planar patches results in a planar graph. 

Let $\P$ be a clean patch with rim $f$. If $(\P,\Omega_f)$ contains a cross, then the union of $H$ and such a cross has a $K_5$-minor, by \sref{addcross}, with every branch set meeting $V(H)$; so that $G_0$ has a $K_6$-minor and \sref{girth5-main} follows. Assume then that $(\P,\Omega_f)$ has no cross and is nonplanar.
Then, $\P = \P_1 \cup \P_2$, $\P_1 \cap \P_2 = \P[D]$ and $|D| \leq 3$ such that
$V(f) \subseteq V(\P_1)$ and $|V(\P_2) \sm V(\P_1)| \geq 2$, by \sref{cross}. Hence, $\{z_0\}\cup D$ is a $k$-disconnector of $G_0$ with $k \leq 4$; contradicting (\sref{girth5-main}.D). It follows that $\P$ is planar so that (\sref{girth5-main}.E) follows.

If $x$ is a vertex-breaker, then let $C$ be the vertices of $H$ cofacial with $x$. $4$-connectivity of $G_1$ implies that every vertex in $H'-\{x\}-C$ is at least $4$-valent in $H'-x$. As $x$ is $3$-valent in this case, by (\sref{ess4}.3), we have that $H'-x$ is a $2$-connected planar graph of girth $\geq 5$ has an embedding in the plane with each vertex not in $X_{H'-x}$ at least $4$-valent, and each vertex in $X_{H'-x}$ at least $3$-valent except for at most $3$ vertices which are at least $2$-valent. By \sref{dis1}, $H'-x$ is does not exist; contradiction.\QED


\end{document}